**From a Student's Question to the Half-Derivative:**
**A Constructivist Teaching Model — Part I**

Davit Kapanadze
Doctor of Pedagogical Sciences
Georgian National University SEU, Tbilisi, Georgia
Associate Member (AMIMA), Institute of Mathematics and its Applications, United Kingdom
d.kapanadze@seu.edu.ge | ORCID: 0009-0006-8112-5084

## Abstract

This paper presents a conceptual, cognitive-constructivist model for teaching the half-derivative. The model starts from a student's question and frames it as the problem of generalizing the familiar operation of differentiation. On a suitably chosen space of monomials, under the assumptions of linearity and of lowering the exponent by one half, the requirement that the operator applied twice give the ordinary derivative produces a recurrence for the coefficients. The freedom of normalization becomes visible, and an additional convention yields a formula compatible with the monomial rule of the Riemann–Liouville half-derivative. The domain of the operator and its relation to the Riemann–Liouville and Caputo operators are treated separately. The mathematical results rest on the author's earlier preprints; the contribution of this paper is their transformation into an instructional pathway organized in phases, with diagnostic questions, an analysis of typical errors, and comparative problems. The model is a design proposal: its effectiveness has not yet been established by systematic empirical research. The interactive Socratic implementation and the detailed diagnostics are the subject of Part II.

**Keywords:** half-derivative; fractional differentiation; constructivist teaching; monomial space; composition of operators; recurrence; normalization; APOS; instructional diagnostics.

### Notation

| Symbol | Meaning |
|---|---|
| $D, D^2$ | ordinary differentiation and its twofold application, $D^2 = D \circ D$ |
| $H$ | the sought half-order operator during the construction (requirement $H \circ H = D$) |
| $H_K = D_K^{1/2}$ | the constructed half-derivative after normalization; the index $K$ distinguishes it from the classical operators |
| $c(\beta)$, $c_0 = c(0)$ | coefficient of $H$ on the monomial $x^\beta$ and the free normalization parameter; after normalization $c_0 = 1/\sqrt{\pi}$ |
| $S, S^+$ | the exponent set $\{0, \frac{1}{2}, 1, \frac{3}{2}, \dots\}$ and its extension $S^+ = S \cup \{-\frac{1}{2}\}$ |
| $V_S, V_S^+$ | spaces of finite linear combinations of the corresponding monomials (original and extended) |
| $D_{RL}^{1/2} = D_{0+}^{1/2}$ | left-sided Riemann–Liouville half-derivative with lower terminal 0 |
| $D_C^{1/2}$ | Caputo half-derivative with lower terminal 0 |

| Symbol | Meaning |
|---|---|
| $D^{\alpha}_{0+}$, $D^{\alpha}_{C}$ | Riemann–Liouville and Caputo derivatives of general order $\alpha > 0$ |
| $J^{\alpha}_{0+}$ | Riemann–Liouville fractional integral of order $\alpha$; $J$ — the ordinary integral $\int_0^x$ |
| $J_K$, $d(\beta)$ | algebraically constructed half-integral and its coefficient (§4.9) |
| $I_{V_S}$, $I_{V_S^+}$ | identity operator on $V_S$ and on $V_S^+$, respectively (§4.9) |
| $\kappa$ | constant term of the solution in Problem C (§6.3) |
| $\Gamma$, $B$ | gamma and beta functions |
| $x$, $t$ | independent variable ($x > 0$ in the construction, $t > 0$ in the classical definitions and the problems) |

# 1. Introduction

## 1.1. From a student's question to a mathematical problem

Knowing how to apply a mathematical rule is not the same as understanding its origin. A student may be able to use a formula and still find it hard to explain where it comes from, which assumptions it rests on, and when it may be applied. The starting point of this paper is that a new concept should not be delivered only as a finished definition; the student should be able to understand the motive for its construction and the steps that lead to it.

*Teaching context.* According to the author's account, in February 2026 three second-year students of the Faculty of Engineering of Georgian National University SEU discussed among themselves whether it makes sense to speak of a "derivative of order one-half". One of them then put the question to the lecturer. The students had not studied fractional calculus and did not know the Riemann–Liouville or Caputo definitions. The episode is given as the pedagogical motivation of the construction, as reconstructed by the author; it is not systematically collected research data, and the students' words are not treated as direct quotations.

The student knows that repeated application of the ordinary derivative gives the second and higher derivatives. This raises a working question: can one construct an operator whose application twice in succession gives ordinary differentiation? Denote the ordinary derivative by $D$ and the sought operator by $H$. The defining requirement is

$$H \circ H = D. \tag{1.1}$$

Equation (1.1) is a problem, not a statement of existence or uniqueness. To solve it, one must choose a class of functions, a possible form of the operator's action, and the conditions under which the composition is admissible. In particular, (1.1) by itself implies neither the linearity of $H$ nor its action on monomials. These conditions are introduced below explicitly as working assumptions.

### 1.2. Mathematical scope of the construction

Wherever monomials with non-integer or negative exponents occur, we work on $x > 0$. Take the set of non-negative integer and half-integer exponents $S = \{0, \frac{1}{2}, 1, \frac{3}{2}, \dots\}$ and the space $V_S = \text{span}\{x^\beta : \beta \in S\}$, where only finite linear combinations are meant. On this space we seek a linear operator that maps the monomial $x^\beta$ to $c(\beta)x^{\beta-1/2}$. Lowering the exponent by one half and the monomial form are working assumptions, not consequences of (1.1).

$$H(x^\beta) = c(\beta)x^{\beta-1/2}, \quad \beta \in S. \qquad (1.2)$$

For $\beta \geq \frac{1}{2}$, applying the operator twice and comparing with the ordinary derivative gives the recurrence

$$c(\beta)c\left(\beta - \frac{1}{2}\right) = \beta, \quad \beta \in S, \beta \geq \frac{1}{2}. \qquad (1.3)$$

The recurrence does not determine the coefficients uniquely: a free normalization parameter $c(0)$ remains. The choice $c(0) = \frac{1}{\sqrt{\pi}}$, which makes the operator agree with the Riemann–Liouville operator on constants, gives the well-known monomial formula

$$H(x^\beta) = \frac{\Gamma(\beta+1)}{\Gamma\left(\beta + \frac{1}{2}\right)} \cdot x^{\beta-1/2}, \quad \beta \in S. \quad (1.4)$$

Since $H(1) = \frac{x^{-1/2}}{\sqrt{\pi}}$ does not belong to $V_S$, the composition can be fully defined only on the extension $V_S^+ = V_S \oplus \text{span}\{x^{-1/2}\}$. On the extended space we additionally set $H(x^{-1/2}) = 0$. Thus $H: V_S^+ \to V_S^+$ and $(H \circ H)f = Df$ for every $f \in V_S$. This identity is not established for all elements of the extended space $V_S^+$: for example, $(H \circ H)(x^{-1/2}) = 0$, whereas $D(x^{-1/2}) = -\frac{1}{2}x^{-3/2}$. Here the difference between the domain of the operator and the domain of validity of the composition identity becomes visible.

Formula (1.4) is not the starting point of the instructional sequence; it enters only after the student has obtained (1.3) step by step, discovered the free parameter, and realized that an additional normalization is needed. The rigorous justification and the comparison with the classical operators are given in the following sections.

### 1.3. Aim and contribution

The aim of the paper is to turn an already known monomial construction into an instructional sequence in which the need for each new step arises from the result or the limitation of the previous one. The mathematical results discussed in the author's earlier preprints [1, 2] are not presented here as new theorems. The contribution of this article is their didactic organization: the student first makes the defining requirement precise, then chooses the assumptions, derives the recurrence, discovers the non-uniqueness of the normalization, links the choice to the classical theory, and checks the limits within which the rule is valid.

The instructional model combines three interrelated components: (a) a stepwise construction of mathematical actions; (b) the specification of anticipated errors and diagnostic questions at each stage; (c) problems comparing the constructed rule with the Riemann–Liouville and Caputo operators. Starting from a student's question is not in itself a new pedagogical device; a similar motivation is used, for example, by Miller and Sugden [3]. The focus here is on a concrete instructional pathway, supported by diagnostics, that leads from composition to the recurrence and from the free normalization to the gamma-function formula. A complete systematic review of existing research was not carried out within this paper.

### 1.4. Instructional logic

The proposed pathway is organized as follows:

*activation of prior knowledge → formalization of the question → choice of assumptions → derivation of the recurrence → computation of the coefficients → discovery of the freedom of normalization → comparison with the classical operators → applied problems.*

Each stage corresponds to an action the student can perform and to a question with which the lecturer checks understanding. For example, when the assumptions are chosen, the student is asked which condition produced the change of exponent; at the normalization stage, why the recurrence alone is not enough; and regarding the domain, when the twofold application of $H$ is legitimately defined.

Thus the theoretical framework (Piaget, Ausubel, Sweller, Vygotsky, APOS, and concept image) is used not as evidence of effectiveness but as a basis for designing concrete instructional actions. Its direct link to these actions is specified in Section 2.

### 1.5. Status and structure of the study

The article is a conceptual-methodological paper. It presents an instructional design proposal and contains no validated experimental results or causal evaluation of effectiveness. The author's informal practical observations and discussions with colleagues were used to refine the idea, not as empirical evidence for the conclusions.

Section 2 connects cognitive-constructivist theories with instructional actions; Section 3 explains the relation to the sources; Section 4 describes the construction and the logical status of the formulas; Sections 5 and 6 present the comparison and the problems; Sections 7–9 are devoted to the proposed classroom implementation, the supporting role of AI, and the diagnostic framework; Section 10 contains the discussion and the limitations. The operation of the interactive program "Sokrate AI — Kapanadze", its support logic, and the detailed diagnostics are discussed in Part II of the study and are not an empirical result of this part.

## 2. Brief theoretical framework

The theories invoked in this section form the basis for designing the instructional model. They are not used as evidence of the model's effectiveness.

### 2.1. Cognitive-constructivist foundations

The model rests on the view that meaningful learning of mathematics is not limited to receiving and memorizing ready-made information. The student should be able to formulate a new problem on the basis of prior knowledge, state a hypothesis, check its consequences, and build a new relationship step by step. From this point of view four theoretical sources are relevant.

**Piaget: assimilation and accommodation.** In Piaget's theory of cognitive development [7], knowledge is the result of the subject's own activity: a new experience either fits into existing schemes (assimilation) or requires their transformation (accommodation). In our case the student first tries to understand the half-order derivative by analogy with familiar operator actions. The student then discovers that the analogy is not sufficient: the class of functions, the action of the operator, and additional conditions must be specified. In this way the initial intuition becomes a formalized problem.

Instructional action (Piaget). The lecturer asks the student to recall the rule for the ordinary derivative and asks: "How should the exponent change if one differentiation is represented by two identical steps?" The student first states a conjecture and then checks it by composition. The diagnostic indicator is whether the student distinguishes a plausible intuition from a proven property.

**Ausubel: meaningful learning.** According to Ausubel's theory [8], new information must be connected to the existing structure of knowledge. Therefore, before instruction begins, it is checked whether the student has the following prior knowledge: computing the derivative of a monomial and understanding the change of exponent; applying an operator repeatedly to a function; comparing coefficients in algebraic expressions; understanding the meaning of a recurrence relation. The construction starts precisely with these familiar actions.

Instructional action (Ausubel). Before starting, the student completes a short four-item diagnostic: (a) compute the derivative of a monomial; (b) explain the composition of two operators; (c) compare the coefficients of two monomials; (d) explain the meaning of a recurrence. Depending on the error, the lecturer returns the student only to the prerequisite step that is needed.

**Sweller: cognitive load.** Sweller's research [9] points to the limited capacity of working memory. If the construction were presented all at once (recurrence, double factorials, normalization and the gamma function together), it would be too demanding for the student. The instructional process is therefore divided into phases: the defining condition for the operator, the possible form of its action, the recurrence relation, specific coefficients, the general pattern, and finally the normalization. In each phase attention is directed to one essential issue.

Instructional action (Sweller). Each phase poses one main task: first the change of exponent, then the product of two coefficients, then the computation of several terms, and finally the normalization. The gamma function is not introduced until the student has a rule for the coefficients discovered through his or her own computations. Answers to the intermediate questions show at which stage the load became excessive.

**Vygotsky: assisted learning.** In Vygotsky's theory [10], what matters are the capabilities that emerge for the student with appropriate support. In the model this is realized through a system of guiding questions. The lecturer does not hand over a ready formula but poses tasks that allow the student to take a new step independently using known actions. If the step cannot be taken, instead of giving the answer the lecturer returns to the prior knowledge causing the difficulty and asks a simpler question. For example, if the student finds it hard to determine by how much $H$ should lower the exponent, one may ask: "If the sum of two equal changes is one unit, how large is each?"

Instructional action (Vygotsky). If the student cannot find the half-unit shift, the lecturer does not give the answer: the lecturer first asks the student to write down the sum of two equal changes and then to find a single change. If the obstacle persists, a simple numerical analogue is offered. Steps performed independently and steps performed with help are recorded separately; this record is not yet a validated measure.

Thus four interrelated principles are distinguished in the model: activation of prior knowledge, active mathematical work by the student, stepwise introduction of new material, and targeted help when needed. The student's question is the starting point of the instructional design here; we do not generalize from one such question to the reactions of all students.

### 2.2. APOS: Action – Process – Object – Schema

The APOS framework of Dubinsky and McDonald [11] describes four mental constructions involved in understanding a mathematical concept: Action, Process, Object, and Schema. In the case of the half-derivative, the framework allows us to specify in advance not only what the student should compute but also how the form of their understanding may change.

Examples of observable APOS behaviour. At the action level the student computes a single case such as $H(x)$; at the process level the student explains why the composition gives $c(\beta)c\left(\beta-\frac{1}{2}\right)=\beta$; at the object level the student reasons about why a choice of normalization exists and which domain is needed for $H^2$; at the schema level the student checks when the normalized rule agrees with the Riemann–Liouville and Caputo rules. This describes observable responses, not a direct measurement of the student's internal state.

**Action.** The student works with individual monomials: recalls the rule for the derivative, considers the composition and, once the rule has been obtained, computes the half-derivative of specific monomials.

**Process.** Beyond individual computations the student sees a general relationship: in the chosen model the exponent decreases by $1/2$, and the coefficients are linked by the recurrence relation (1.3).

**Object.** The half-derivative is no longer only a rule for transforming a monomial. It is an operator in its own right, with a domain, conditions of validity, and relations to other operators; it can be applied repeatedly, inverted, and compared.

**Schema.** The student combines in one system the monomial rule, the recurrence, the normalization, the gamma-function formula, the domains, and the comparison with the classical operators; on this basis the student can use the operator in problems and assess the legitimacy of the result (Table 1).

*Table 1. APOS stages in the proposed teaching of the half-derivative*

| Stage | Expected student activity | Diagnostic question |
|---|---|---|

| Stage | Expected student activity | Diagnostic question |
|---|---|---|
| Action | computes $H(x)$ and $H(x^2)$ by the given rule | What does the operator do to a specific monomial? |
| Process | combines the change of exponent and coefficient into a general pattern | What common rule underlies these examples? |
| Object | investigates $H \circ H$, the normalization, and the domain | What must be specified for $H$ to be treated as an operator in its own right? |
| Schema | combines the constructed and classical operators, the problems, and the conditions for composition | When do the actions of the operators coincide and when do they differ? |

**From process to object.** In APOS terms the hardest step is usually encapsulation: a process the student can carry out becomes an object on which new actions are performed. Three situations are designed in the sequence to support this step.

1. *The free parameter $c_0$ (Section 4.4).* The composition requirement is satisfied not by a single rule but by a whole family of rules, one for each $c_0 \neq 0$. The question "What does the rule do to $x^\beta$?" is replaced by the question "Which operator do we choose from the family, and on what grounds?" The operator itself becomes the object of reasoning.
2. *The domain of the composition (Sections 4.2, 4.7 and 4.8).* The question whether $H \circ H$ is defined on a given function treats $H$ as an object with a domain and a range. This is where the main observation of §1.2 arises: the domain of the operator ($V_S^+$) and the domain of validity of the identity ($V_S$) are different.
3. *The operator in an equation (Problems B and C of Section 6).* When the function is unknown and the operator is given, $H$ is used as an object of the problem in its own right. What matters is not only the formal manipulation but also checking whether inversion is admissible, what the kernel of the operator is, and whether the solution may be non-unique.

The internal mental construction of a particular student cannot be established from a correct answer alone. A plausible observable indicator of the transition is a change in the type of the student's statements: from statements about results ("$H(x^2)$ equals …") to statements about the operator ("$H$ is not uniquely determined by the composition condition", "$H$ agrees with the Caputo operator on monomials with positive exponent but not on constants"). These indicators are hypothetical; their reliability requires empirical testing. The rubric in Section 9 is designed precisely to record such statements.

### 2.3. Concept image and its limits

Tall and Vinner's notion of concept image [12] allows us to distinguish the student's first intuition — "half a step of the ordinary derivative" — from the rigorous mathematical object. The intuition is a good starting point for raising the question, but it is not a definition and may even be misleading: $D^{1/2}$ is not obtained by substituting the number $1/2$ into the formula for the derivative, and $H \circ H = D$ does not automatically determine a unique operator on all functions.

Instructional action (concept image). The lecturer first records the student's explanation "half a step" and then asks the student to extract from this phrase only the hypothesis and to list what it does not establish: the existence of the operator, its uniqueness, its linearity, and its domain. At the end of the lesson the student compares their initial explanation with the definition obtained.

The intuitive picture should therefore be gradually refined and linked to formal requirements. The student should realize that constructing the operator requires specifying the class of functions on which it acts, its algebraic form, the domain of validity of the composition condition, the normalization, and the relation to the classical fractional derivatives. In this way the initial concept image is enriched with new connections and limitations.

Orton's study [13] showed that many difficulties in learning the ordinary derivative are conceptual rather than purely procedural. In the case of the half-derivative this gap between rule and meaning may widen further, because a familiar concept is extended to a case in which new conditions of validity and new limitations arise. The lecturer's task is not to reject the intuition but to use it as an initial hypothesis that is subsequently tested and refined.

## 3. Type of study and relation to previous work

### 3.1. Type of study: a design proposal

The paper is a conceptual-methodological study. Its outcome is an instructional design proposal based on a comparative analysis of mathematically verified examples. The design consists of four components, each described so that it can be tested separately in later empirical research.

1. **Instructional phases** (Section 4). For each phase the prior knowledge, the expected student action, the possible obstacle, and the diagnostic question are specified.
2. **Logical status of mathematical statements.** Every formula the student meets is labelled according to its logical status: *requirement* (e.g., $H \circ H = D$, (1.1)); *assumption* (monomial form and linearity, (4.3)); *derived result* (the recurrence (1.3)); *convention* ($c_0 = 1/\sqrt{\pi}$); *compact notation* (the gamma formula (1.4)). This labelling also serves an instructional purpose: at every step the student should know whether a formula has been derived or chosen.
3. **Problems** (Section 6). Each problem is solved both classically and constructively; its role in the sequence is stated separately.
4. **Assessment criteria** (Section 9). The rubric describes individual criteria. Until it has been validated, results should be reported as a criterion profile rather than as a total score.

These components have not been tested as an instructional model, and we do not present them as such.

### 3.2. Mathematical basis and the didactic contribution

The mathematical basis of the paper is the author's two preprints [1, 2]; these results are not presented here as rediscovered theorems. The theoretical framework relies on the sources discussed in Section 2 [7–13]. Table 2 separates the reused mathematical content from the didactic contribution of this paper.

*Table 2. Separation of reused mathematical content and the didactic contribution of this paper*

| **Taken from [1, 2] (mathematics)** | **Contribution of this paper (didactics)** |
|---|---|
| the recurrence relation $c(\beta)c\left(\beta-\frac{1}{2}\right)=\beta$ on the graded monomial space | a task divided into phases in which the student derives the recurrence from $H \circ H = D$ on his or her own, with diagnostic questions for each phase |
| two families of coefficients, double factorials, and the free parameter $c_0$ | the free parameter as a planned didactic obstacle that motivates normalization and the transition from process to object |
| the normalization $c_0 = 1/\sqrt{\pi}$ and the gamma-function form | normalization as an explicitly stated convention, separated from the results derived from the recurrence |
| the monomial rule of general order and the conditions for composition | a comparison of three operators on the same set of monomials, arranged so that both agreement and difference are visible |
| application of the construction to specific functions | seven problems with parallel classical and constructive solutions, each with a stated instructional role |
| — | a catalogue of anticipated errors with the lecturer's responses; brief principles for the use of AI; assessment criteria |

The individual elements of the right-hand column are not new in themselves: diagnostic questions, an error catalogue, and a rubric are standard tools of instructional design. What is new concerns the three points listed in §1.3 — the route from composition to the recurrence, the discovery of the freedom of normalization, and the organization of these steps into an instructional sequence — that is, how these tools are selected and linked for one specific, non-traditional mathematical construction.

### 3.3. Relation to the usual order of presenting fractional calculus

In standard textbooks [14–18] the subject starts with definitions: Cauchy's formula for repeated integration is generalized to the Riemann–Liouville fractional integral, the Riemann–Liouville and Caputo derivatives are defined through it, and the monomial formula is obtained from properties of the beta and gamma functions. This route is mathematically sound and necessary for the general theory, but from the very beginning it requires a fairly advanced knowledge of integral calculus and special functions.

Historically the monomial rule is older than these definitions: Lacroix obtained the half-derivative of a power by replacing the factorials in the ordinary power rule with gamma functions [19]. This formal substitution gives the appropriate coefficient, but it does not answer the instructional question: how can a student arrive at this coefficient through familiar actions?

The proposed sequence takes a third position. Like Lacroix's, it starts from monomials, but it does not postulate the gamma coefficient: the coefficient is derived from the composition requirement, the non-uniqueness becomes visible, and the classical theory enters exactly at the moment when a choice has to be made — when the student is already able to ask which operator to choose from the family and why. The differences are summarized in Table 3.

*Table 3. Definition-first presentation and the proposed constructive sequence*

| **Aspect** | **Definition-first presentation [14–18]** | **Proposed constructive sequence** |
|---|---|---|
| Starting point | the fractional integral as a generalization of repeated integration | the student's question and the requirement $H \circ H = D$ |
| Prior knowledge | improper integrals, beta and gamma functions | the power rule, composition of operators, comparison of coefficients |
| Origin of the monomial formula | from the integral definition | from the chosen monomial form and the composition requirement |
| Role of the gamma function | a starting tool, built into the definition | at the end, as a compact notation for coefficients already computed |
| Normalization | built into the definition and not visible as a choice | an explicit convention after the free parameter has been discovered |
| Riemann–Liouville and Caputo | two definitions that are compared afterwards | the outcome of the question of which action on constants should be reproduced |
| Class of functions | broad (locally integrable or absolutely continuous functions) | narrow (monomials and their finite linear combinations) |
| Main limitation | high entry threshold; the origin of the formula may remain opaque | does not by itself reach functions outside the monomial class |

The two approaches complement rather than compete with each other. The merit of the classical presentation is its direct connection with the general theory; the possible didactic merit of the constructive sequence is that its first steps are carried out with already familiar algebraic and differential operations. The initial simplicity of the method does not, however, mean that the whole theory is simplified: the monomial construction is defined on a restricted class, and passing to general functions requires further analysis. The sequence is therefore regarded as a conceptual introduction that precedes the integral definitions and motivates them, not as a replacement for them. The relation to the work of Miller and Sugden [3] is discussed in §1.3; a systematic review of the didactic literature on fractional calculus is beyond the scope of this paper.

**4. The constructive pathway: from the student's question to the gamma formula**

This section describes the instructional phases (Table 4). For each phase the prior knowledge, the expected student action, the possible obstacle, and the diagnostic question are given. The main formulas carry their logical status in square brackets (§3.1): [**requirement**], [**assumption**], [**derived result**], [**convention**] or [**compact notation**]. The lecturer introduces the topic only in general terms at the outset; the steps in the construction are small — each phase differs from the previous one by a single new idea.

*Table 4. Instructional phases of the construction*

| **Phase** | **Student action** | **Possible obstacle** | **Diagnostic question** |
|---|---|---|---|

| Phase | Student action | Possible obstacle | Diagnostic question |
|---|---|---|---|
| 1. Formalizing the question (§4.1) | formulates the requirement $H \circ H = D$ | takes the equation as a statement of existence | What are we looking for, and what do we not yet know? |
| 2. Choosing the form (§4.2) | chooses the monomial form, finds the shift $1/2$ | treats an assumption as a derived result | Does this follow from $H \circ H = D$, or did we choose it? |
| 3. Recurrence (§4.3) | obtains $c(\beta)c\left(\beta - \frac{1}{2}\right) = \beta$ | forgets the condition of the derivation $\beta \geq \frac{1}{2}$ | For which $\beta$ is this derivation legitimate? |
| 4. Families and $c_0$ (§4.4) | computes the coefficients, sees two families | expects a unique answer | Why is the recurrence not violated when $c_0$ is changed? |
| 5. Normalization (§4.5) | chooses $c_0 = 1/\sqrt{\pi}$ | takes the normalization to be a consequence of the recurrence | Which additional condition did we introduce, and why this one? |
| 6. Gamma form and domain (§4.6–4.7) | writes the general formula, fills in the table, checks the domain | assumes the identity also holds on the extended space | Where is $H_K$ defined, and where does $H_K \circ H_K = D$ hold? |
| 7. Extension (§4.8–4.9; for interested students) | general orders, the half-integral | applies the composition rule unconditionally | Why is adding the exponents not enough? |

### 4.1. The question and the operator hypothesis

The lecturer starts from a familiar action, $D(Df) = D^2 f$, and poses the following question: can one differentiation be represented as the result of two identical operator steps? This leads to the defining condition

$$H \circ H = D \qquad (4.1)$$

[**requirement**], where $H$ is the sought half-order operator. As noted in §1.1, (4.1) does not assert that such an operator exists or is unique; it is a problem, and to solve it the class of functions, the rule of action, and additional conditions must be specified. The main question of the phase is: which functions should we start from so that the results can be checked?

A note on notation: during the construction we write $H$; the operator obtained after the normalization and the domain have been fixed is denoted by $H_K = D_K^{1/2}$ (§1.2), to distinguish it from the Riemann–Liouville (RL) and Caputo (C) operators.

### 4.2. The monomial space and the starting point

As the initial class we choose the monomials $x^\beta$, $x > 0$, where

$$S = \left\{0,\ \frac{1}{2},\ 1,\ \frac{3}{2},\ 2,\ \dots\right\}, \qquad (4.2)$$

and their finite linear combinations form the space $V_S$ of §1.2. The choice also has a methodological basis: the student knows $D(x^\beta) = \beta x^{\beta-1}$, so the construction of the new operator starts with familiar algebra.

The lecturer asks: "If the ordinary derivative lowers the exponent by one, by how much must an operator lower it if applying it twice equals one differentiation?" To discuss this question we explicitly restrict the class of operators under consideration: $H$ is linear, maps monomials to monomials, and lowers the exponent each time by the same amount $s$. These properties do not follow from (4.1); they are assumptions of the construction declared in advance. Applied twice, the exponent decreases by $2s$; comparison with differentiation gives $2s = 1$, i.e. $s = \frac{1}{2}$. Thus

$$H(x^\beta) = c(\beta)\, x^{\beta-1/2}, \qquad \beta \in S, \qquad (4.3)$$

[**assumption**], where $c(\beta)$ is a coefficient still to be found; finding it is the task of the next phase. On finite sums $H$ is extended linearly. We note right away that the first application to a constant gives $x^{-1/2}$, which lies outside $S$; we return to this point in §4.7.

### 4.3. Independent derivation of the recurrence relation

The student receives the task: "Apply $H$ twice and compare the result with the ordinary derivative." If $\beta \in S$ and $\beta \geq \frac{1}{2}$, the intermediate exponent $\beta - \frac{1}{2}$ again belongs to $S$, and applying (4.3) twice gives

$$H\left(H(x^\beta)\right) = c(\beta)\, c\left(\beta - \frac{1}{2}\right) x^{\beta-1}. \qquad (4.4)$$

On the other hand, $D(x^\beta) = \beta x^{\beta-1}$. The requirement (4.1) means equality of the coefficients:

$$c(\beta)\, c\left(\beta - \frac{1}{2}\right) = \beta, \qquad \beta \in S,\ \beta \geq \frac{1}{2} \qquad (4.5)$$

[**derived result**]. The relation is obtained through three familiar actions: composition of operators, comparison of exponents, and comparison of coefficients. This is the central result of the construction.

The condition of the derivation matters. For $\beta = 0$ the first application of $H$ leaves the initial class, so the composition on a constant is not covered by (4.5) and is treated separately (§4.7). From an instructional point of view the student sees here that writing down a symbolic equality does not yet make its use legitimate on all functions.

### 4.4. Two families and the free normalization

Let $c(0) = c_0$. From (4.5) at $\beta = \frac{1}{2}$: $c\left(\frac{1}{2}\right)\, c_0 = \frac{1}{2}$, that is,

$$c\left(\frac{1}{2}\right) = \frac{1}{2c_0}, \qquad c_0 \neq 0. \qquad (4.6)$$

Thus $c_0 \neq 0$ is not an additional assumption but a necessary consequence of the recurrence. Then for $\beta = 1, \frac{3}{2}, 2, \dots$ one obtains step by step:

$$c(0) = c_0, \quad c\left(\frac{1}{2}\right) = \frac{1}{2c_0}, \quad c(1) = 2c_0, \quad c\left(\frac{3}{2}\right) = \frac{3}{4c_0}, \quad c(2) = \frac{8c_0}{3},$$
$$c\left(\frac{5}{2}\right) = \frac{15}{16c_0}, \quad c(3) = \frac{16c_0}{5}. \tag{4.7}$$

In this sequence the student notices two families: for integer exponents $c_0$ is in the numerator, for half-integer exponents in the denominator. The pattern is written with double factorials:

$$c(n) = c_0 \frac{(2n)!!}{(2n-1)!!}, \qquad c\left(n + \frac{1}{2}\right) = \frac{(2n+1)!!}{2c_0\,(2n)!!}, \qquad n = 0,1,2,\dots, \tag{4.8}$$

where $0!! = (-1)!! = 1$. The double factorial does not enter here as a ready-made formula — it records a pattern discovered in coefficients already computed.

The central diagnostic question of the phase is: "If we change $c_0$, why does (4.5) remain valid?" Answer: in neighbouring coefficients $c_0$ appears once in the numerator and once in the denominator, and cancels in the product. In this way the student discovers that the composition requirement determines a *family* of operators rather than a single operator. In APOS terms this is the first situation of the transition from process to object (§2.2): the reasoning shifts from computing coefficients to choosing an operator. The Wallis product shows an asymptotic relation between these families [1], but by itself does not determine the value of $c_0$.

### 4.5. Normalization: why is an additional condition needed?

After the recurrence has been solved a new problem arises: which operator in the family should we choose? At this moment the classical theory enters. The lecturer gives the student a known result: the left-sided Riemann–Liouville half-derivative with lower terminal 0 gives on a constant

$$D_{RL}^{1/2}\, 1 = \frac{1}{\sqrt{\pi x}}, \qquad x > 0. \tag{4.9}$$

For the constructed operator $H(1) = c_0\, x^{-1/2}$. If we require that it agree with the RL operator on constants,

$$H(1) = D_{RL}^{1/2}\, 1, \tag{4.10}$$

we obtain the normalization

$$c_0 = \frac{1}{\sqrt{\pi}} \tag{4.11}$$

**[convention]**. In the scenario (4.11) remains a convention; the derivation of (4.9) from the integral definition is offered to interested students (§4.8). The main point is that the student realizes that (4.11) does not follow from the recurrence; it is chosen by an additional compatibility condition.

At the same moment the difference from the Caputo operator appears. The Caputo half-derivative of a constant is zero, $D_C^{1/2}1 = 0$, whereas (4.6) requires $c_0 \neq 0$. Consequently, the Caputo behaviour on constants does not fit into this family. Illustration: $D_C^{1/2}x^{1/2} = \frac{\sqrt{\pi}}{2}$, and the Caputo half-derivative of this constant is zero, while the ordinary derivative is $\frac{1}{2\sqrt{x}}$. Hence applying the Caputo half-derivative twice to $x^{1/2}$ does not satisfy (4.1). As the outcome of the phase, the student clearly separates two things: what follows from the composition requirement and what is determined by the additionally chosen normalization.

### 4.6. The main formula and its gamma-function form

Substituting (4.11) into (4.8), both families take a concrete form:

$$c(n) = \frac{(2n)!!}{(2n-1)!!\,\sqrt{\pi}}, \qquad c\left(n+\frac{1}{2}\right) = \frac{(2n+1)!!\,\sqrt{\pi}}{2\,(2n)!!}.$$

By the properties $\Gamma(z+1) = z\Gamma(z)$ and $\Gamma\left(\frac{1}{2}\right) = \sqrt{\pi}$, the two families merge into one formula:

$$H_K\left(x^\beta\right) = \frac{\Gamma(\beta+1)}{\Gamma\left(\beta+\frac{1}{2}\right)}\,x^{\beta-1/2}, \qquad \beta \in S,\ x > 0 \qquad (4.12)$$

[**compact notation**]. This is formula (1.4) of §1.2, and it coincides with the monomial rule of the left-sided RL half-derivative with lower terminal 0. Thus the construction does not produce a new classical derivative; its instructional value lies in the route by which the student goes from familiar algebra to a known formula:

> **composition → choice of the shift → recurrence → two families → free parameter → normalization → gamma formula**

### 4.7. Table of values and the domain of the operator

Formula (4.12) allows us to compute the half-derivative on the monomials of $S$ (Table 5). The table is not only a display of results: the lecturer invites the student to fill in several cells, to check the product of neighbouring coefficients against (4.5), and to compare the result with (4.12).

*Table 5. Half-derivatives of monomials before and after normalization*

| $\beta$ | $c(\beta)$, **with free** $c_0$ | $c(\beta)$, $c_0 = 1/\sqrt{\pi}$ | $H_K\left(x^\beta\right)$ |
|---|---|---|---|
| $0$ | $c_0$ | $\frac{1}{\sqrt{\pi}}$ | $\frac{1}{\sqrt{\pi x}}$ |
| $\frac{1}{2}$ | $\frac{1}{2c_0}$ | $\frac{\sqrt{\pi}}{2}$ | $\frac{\sqrt{\pi}}{2}$ |
| $1$ | $2c_0$ | $\frac{2}{\sqrt{\pi}}$ | $\frac{2\sqrt{x}}{\sqrt{\pi}}$ |

| $\beta$ | $c(\beta)$, **with free** $c_0$ | $c(\beta)$, $c_0 = 1/\sqrt{\pi}$ | $H_K(x^\beta)$ |
|---|---|---|---|
| $\frac{3}{2}$ | $\frac{3}{4c_0}$ | $\frac{3\sqrt{\pi}}{4}$ | $\frac{3\sqrt{\pi}}{4}x$ |
| $2$ | $\frac{8c_0}{3}$ | $\frac{8}{3\sqrt{\pi}}$ | $\frac{8\,x^{3/2}}{3\sqrt{\pi}}$ |
| $\frac{5}{2}$ | $\frac{15}{16c_0}$ | $\frac{15\sqrt{\pi}}{16}$ | $\frac{15\sqrt{\pi}}{16}x^2$ |
| $3$ | $\frac{16c_0}{5}$ | $\frac{16}{5\sqrt{\pi}}$ | $\frac{16\,x^{5/2}}{5\sqrt{\pi}}$ |

For example, the student checks: $H_K(x) = \frac{2\sqrt{x}}{\sqrt{\pi}}$ and $H_K\big(H_K(x)\big) = \frac{2}{\sqrt{\pi}} \cdot \frac{\sqrt{\pi}}{2} = 1 = D(x)$. In this way (4.1) is verified on an individual monomial.

**Extending the domain.** On a constant $H_K(1) = \pi^{-1/2}x^{-1/2}$, and the exponent $-\frac{1}{2}$ does not belong to $S$. To make sense of the second application we extend the exponent set: $S^+ = S \cup \{-\frac{1}{2}\}$. On a constant the composition requirement is $H\big(H(1)\big) = D(1) = 0$; since $c_0 \neq 0$, we must have $c\left(-\frac{1}{2}\right) = 0$, that is,

$$H_K\big(x^{-1/2}\big) = 0. \qquad (4.13)$$

This value does not depend on the normalization and is also the natural continuation of (4.12), since $1/\Gamma(0) = 0$.

Domains of the operator and of the identity. Defining $H_K$ on the extended space does not mean that $H_K^2 = D$ holds for every element of the extended space. The identity is established only on the original space $V_S$; the difference is made explicit by the check below.

Let $S = \{0, \frac{1}{2}, 1, \frac{3}{2}, \dots\}$ and $V_S = \mathrm{span}\{x^\beta : \beta \in S\}$, where $x > 0$ and only finite linear combinations are considered. Introduce the extended space $V_S^+ = V_S \oplus \mathrm{span}\{x^{-1/2}\}$.

We define the normalized operator on $V_S^+$: $H_K(x^\beta) = \frac{\Gamma(\beta+1)}{\Gamma\left(\beta+\frac{1}{2}\right)}x^{\beta-1/2}$ $(\beta \in S)$, and $H_K(x^{-1/2}) = 0$. By linear extension $H_K(V_S^+) \subseteq V_S^+$; hence the composition $H_K^2$ is defined on $V_S^+$.

For every $f \in V_S$ we have $H_K^2 f = Df$. In particular, $H_K(1) = \frac{x^{-1/2}}{\sqrt{\pi}}$ and $H_K^2(1) = 0 = D(1)$. But for the additional monomial in the extended space, $H_K^2(x^{-1/2}) = 0$, whereas $D(x^{-1/2}) = -\frac{1}{2}x^{-3/2} \neq 0$.

Thus the domain of $H_K$ is $V_S^+$, while the domain of validity of the identity $H_K^2 = D$ proved here is $V_S$. This distinction is essential for understanding the operator as a mathematical object.

Diagnostic question: why is it necessary to add $x^{-1/2}$ to the domain, and why can we not conclude afterwards that $H_K^2 = D$ on every element of the extended space?

### 4.8. Towards more general orders

Once the half-derivative has been constructed, a natural question arises: can one differentiation be split into three or more equal steps? This leads the student to $\frac{1}{3}, \frac{2}{3}, \frac{3}{2}$ and more general orders. The functional continuation of the coefficients and the domain of the composition are worked out in [2]; here the general rule is presented as the mathematical horizon of the construction, not as an independent derivation. This subsection is intended for interested students.

**The general monomial rule.** For the RL operator with lower terminal 0

$$D_{0+}^{\alpha} x^{\beta} = \frac{\Gamma(\beta+1)}{\Gamma(\beta-\alpha+1)}\, x^{\beta-\alpha}, \qquad \alpha > 0,\ \beta > -1,\ x > 0. \tag{4.14}$$

Poles of the gamma function. When interpreting the coefficient in (4.14), one uses the analytic continuation of $\frac{1}{\Gamma(z)}$ defined on the whole complex plane (which vanishes at z = 0, –1, –2, …); this does not mean that Γ(z) itself is finite there.

**Check:** $\beta = -\frac{1}{2}$. (4.13) can be compared with the classical definition. The half-order integral

$$J_{0+}^{1/2} f(x) = \frac{1}{\sqrt{\pi}} \int_0^x \frac{f(t)}{\sqrt{x-t}}\, dt$$

gives for $f(t) = t^{-1/2}$ the value $J_{0+}^{1/2} t^{-1/2} = \frac{1}{\sqrt{\pi}} \int_0^x \frac{dt}{\sqrt{t(x-t)}} = \sqrt{\pi}$, hence $D_{0+}^{1/2} x^{-1/2} = D\big(\sqrt{\pi}\big) = 0$. The same result is confirmed in three ways: by the composition requirement, by integral computation, and by the appropriate interpretation of the gamma formula. The same integral also yields (4.9): $J_{0+}^{1/2} 1 = 2\sqrt{x/\pi}$, whence $D_{0+}^{1/2} 1 = D\big(2\sqrt{x/\pi}\big) = 1/\sqrt{\pi x}$.

**The composition rule and its limitations.** The formal rule

$$D_{0+}^{\alpha}\big(D_{0+}^{\gamma} f\big) = D_{0+}^{\alpha+\gamma} f \tag{4.15}$$

is not an unconditional law. On monomials, when $\beta > -1$, $\alpha, \gamma > 0$ and the intermediate exponent satisfies $\beta - \gamma > -1$, applying (4.14) twice gives

$$D_{0+}^{\alpha}\big(D_{0+}^{\gamma} x^{\beta}\big) = \frac{\Gamma(\beta+1)}{\Gamma(\beta-\gamma+1)} \cdot \frac{\Gamma(\beta-\gamma+1)}{\Gamma(\beta-\alpha-\gamma+1)}\, x^{\beta-\alpha-\gamma} = D_{0+}^{\alpha+\gamma} x^{\beta},$$

and the cancellation of the gamma functions is legitimate precisely because the intermediate monomial lies in the domain. If $\beta - \gamma \le -1$, the rule may fail. For example, $D_{0+}^{1/2}\big(D_{0+}^{1/2} 1\big) = D_{0+}^{1/2} \frac{1}{\sqrt{\pi x}} = 0 = D(1)$, but $D_{0+}^{1/2}(D\,1) = 0$, whereas $D_{0+}^{3/2} 1 = -\frac{x^{-3/2}}{2\sqrt{\pi}} \neq 0$. Checking the order of the operators and the domain is therefore essential.

Diagnostic question: why is adding the exponents not sufficient to justify the equality $D^{\alpha} D^{\gamma} = D^{\alpha+\gamma}$?

Precise conditions. The general monomial rule $D_{0+}^{\alpha} x^{\beta} = \frac{\Gamma(\beta+1)}{\Gamma(\beta-\alpha+1)} x^{\beta-\alpha}$ is considered under the conditions $x > 0$, $\alpha > 0$ and $\beta > -1$, from the classical Riemann–Liouville definition and with the corresponding pointwise meaning. If $\beta - \alpha + 1$ is a non-positive integer, the coefficient should be read as $\Gamma(\beta+1) \cdot \left[\frac{1}{\Gamma(\beta-\alpha+1)}\right]$.

The composition $D_{0+}^{\alpha}\left(D_{0+}^{\gamma} x^{\beta}\right) = D_{0+}^{\alpha+\gamma} x^{\beta}$ is obtained formally when the intermediate monomial is admissible for applying the classical formula, namely $\beta > -1$ and $\beta - \gamma > -1$, and no exception arising from a zero intermediate coefficient occurs. In this case $(\alpha, \gamma > 0)$ the cancellation of the gamma coefficients is legitimate. Otherwise both sides must be computed directly and their definedness checked separately; merely adding $\alpha + \gamma$ is not a justification of the composition.

### 4.9. The half-integral and conditions for operator inversion on monomial spaces

Problems B and D of Section 6 require the half-order integral. For the sake of a uniform sequence we also construct its monomial action algebraically. We seek an operator $J_K$ that raises the exponent by $\frac{1}{2}$: $J_K\left(x^{\beta}\right) = d(\beta)\, x^{\beta+1/2}$, and require $H\left(J_K\left(x^{\beta}\right)\right) = x^{\beta}$. From the monomial forms $d(\beta)\, c\left(\beta + \frac{1}{2}\right) = 1$, that is $d(\beta) = 1/c\left(\beta + \frac{1}{2}\right)$, and with the normalized coefficients

$$J_K\left(x^{\beta}\right) = \frac{\Gamma(\beta+1)}{\Gamma\left(\beta + \frac{3}{2}\right)}\, x^{\beta+1/2}. \qquad (4.16)$$

This is the monomial rule of the left-sided RL half-integral with lower terminal 0. It can be checked in two ways: $H\left(J_K\left(x^{\beta}\right)\right) = x^{\beta}$, and the composition of two half-integrals

$$J_K\left(J_K\left(x^{\beta}\right)\right) = d(\beta)\, d\left(\beta + \frac{1}{2}\right) x^{\beta+1} = \frac{x^{\beta+1}}{\beta+1} = \int_0^x t^{\beta}\, dt, \qquad \beta > -1.$$

Thus on the monomial class $J_K \circ J_K = J$, where $J$ is the ordinary integral with lower terminal 0.

Conditions for inversion. $H_K$ and the monomial map $J_K$ of the half-integral act in different directions. Equalities between their compositions must be stated with the corresponding spaces indicated; we do not claim invertibility on a general function space.

Recall $V_S$ and $V_S^{+}$ from §4.7. Define $J_K\colon V_S^{+} \to V_S$ by $J_K\left(x^{\beta}\right) = \frac{\Gamma(\beta+1)}{\Gamma\left(\beta+\frac{3}{2}\right)} x^{\beta+1/2}$, $\beta \in S \cup \{-\frac{1}{2}\}$, and extend it linearly. In particular, $J_K(x^{-1/2}) = \sqrt{\pi}$.

Then a direct check on monomials gives $H_K \circ J_K = I_{V_S^{+}}$ and $J_K \circ H_K = I_{V_S}$, where the second composition is considered on $V_S$. The second equality does not extend to $V_S^{+}$: $J_K(H_K(x^{-1/2})) = J_K(0) = 0 \neq x^{-1/2}$.

Consequently, $H_K\colon V_S \to V_S^{+}$ and $J_K\colon V_S^{+} \to V_S$ are mutually inverse bijections. On the extended space, however, $H_K\colon V_S^{+} \to V_S^{+}$ is not injective, since $H_K(x^{-1/2}) = 0$. Accordingly, inversion and uniqueness of solutions on a general space require further analysis.

Diagnostic question: why must the function obtained by solving an operator equation be substituted back into the original equation?

## 5. Classical definitions and a comparison of the three approaches

### 5.1. Riemann–Liouville and Caputo of order one-half

In Section 4 the operator $H_K$ was constructed from the composition requirement and an additional normalization. We now compare it with the classical definitions. The aim is not to prove any superiority of $H_K$, but to establish on which functions the results of the three operators coincide, where they differ, and what role the domain plays. In the classical definitions the variable is denoted by $t$, $t > 0$, and the lower terminal of all operators is $0$.

The Riemann–Liouville half-order integral is

$$J_{0+}^{1/2} f(t) = \frac{1}{\sqrt{\pi}} \int_0^t \frac{f(s)}{\sqrt{t-s}}\, ds. \qquad (5.1)$$

For a monomial $f(t) = t^\beta$ it exists when $\beta > -1$ (integrability near zero). The half-derivatives are defined in two opposite orders:

$$D_{RL}^{1/2} f = \frac{d}{dt}\left[J_{0+}^{1/2} f\right], \qquad D_C^{1/2} f = J_{0+}^{1/2}[f']. \qquad (5.2)$$

Riemann–Liouville (RL) first integrates and then differentiates; Caputo (C) does the opposite. Changing the order of operations leads both to different results and to different requirements on the regularity of the function. A standard sufficient condition for the Caputo definition is absolute continuity of $f$ on $[0,T]$; in the RL case the fractional integral must exist and be differentiable — the precise requirements depend on whether the operator is considered on a function space or pointwise for $t > 0$.

On monomials with positive exponent both classical operators give the same result:

$$D_{RL}^{1/2} t^\beta = D_C^{1/2} t^\beta = \frac{\Gamma(\beta+1)}{\Gamma\left(\beta+\frac{1}{2}\right)}\, t^{\beta-1/2}, \qquad \beta > 0. \qquad (5.3)$$

With the normalization $c_0 = 1/\sqrt{\pi}$, the constructed operator $H_K$ coincides with the RL monomial rule for *all* exponents of $S$, including $\beta = 0$ (4.12). For negative exponents ($-1 < \beta < 0$) the RL monomial expression may be considered pointwise for $t > 0$, although the function-space conditions have to be checked separately; the Caputo definition (5.2) does not apply to these monomials, since $f'$ is not integrable near zero. For a direct comparison of the three operators we therefore take monomials with $\beta \geq 0$ (Table 6).

Precise range of agreement. A direct comparison of all three operators must be carried out on the common part of their domains. For monomials with positive half-integer and integer exponents $\beta \in S$, $\beta > 0$, $x > 0$, the operators $H_K$, $D_{0+}^{1/2}$ (Riemann–Liouville) and $D_C^{1/2}$ (Caputo) give the same result. On a constant $H_K$ agrees with Riemann–Liouville but not with Caputo. The agreement of $H_K$ with the classical monomial rule does not by itself mean that the operator is defined on every function on which the classical derivative exists.

*Table 6. Comparison of the Riemann–Liouville, Caputo and constructed operators on monomials (lower terminal* $0; t > 0$*)*

| $f(t)$ | **Riemann–Liouville** | **Caputo** | **Constructed** $H_K$ |
|---|---|---|---|
| 1 | $\frac{1}{\sqrt{\pi t}}$ | 0 | $\frac{1}{\sqrt{\pi t}}$ |
| $\sqrt{t}$ | $\frac{\sqrt{\pi}}{2}$ | $\frac{\sqrt{\pi}}{2}$ | $\frac{\sqrt{\pi}}{2}$ |
| $t$ | $\frac{2\sqrt{t}}{\sqrt{\pi}}$ | $\frac{2\sqrt{t}}{\sqrt{\pi}}$ | $\frac{2\sqrt{t}}{\sqrt{\pi}}$ |
| $t^2$ | $\frac{8\,t^{3/2}}{3\sqrt{\pi}}$ | $\frac{8\,t^{3/2}}{3\sqrt{\pi}}$ | $\frac{8\,t^{3/2}}{3\sqrt{\pi}}$ |
| $t^3$ | $\frac{16\,t^{5/2}}{5\sqrt{\pi}}$ | $\frac{16\,t^{5/2}}{5\sqrt{\pi}}$ | $\frac{16\,t^{5/2}}{5\sqrt{\pi}}$ |

The didactic function of the table is that the student, besides writing down the results, explains *why* the operators agree in one case and differ in another. Two control examples serve this purpose: one of agreement (§5.2) and one of difference (§5.3).

### 5.2. An example of agreement: $f(t) = t^2$ in three ways

This is the only place in the paper where the half-derivative of $t^2$ is computed in full, in all three ways; Problem A of Section 6 uses the result obtained here and checks only the composition.

**Riemann–Liouville.** Substituting $s = tu$:

$$J_{0+}^{1/2} t^2 = \frac{t^{5/2}}{\sqrt{\pi}} \int_0^1 u^2\,(1-u)^{-1/2}\,du = \frac{t^{5/2}}{\sqrt{\pi}}\,B\left(3,\frac{1}{2}\right) = \frac{16\,t^{5/2}}{15\sqrt{\pi}}, \qquad (5.4)$$

since $B\left(3,\frac{1}{2}\right) = \Gamma(3)\Gamma\left(\frac{1}{2}\right)/\Gamma\left(\frac{7}{2}\right) = \frac{16}{15}$. Differentiating, $D_{RL}^{1/2} t^2 = \frac{d}{dt}\left[\frac{16\,t^{5/2}}{15\sqrt{\pi}}\right] = \frac{8\,t^{3/2}}{3\sqrt{\pi}}$.

**Caputo.** $f'(t) = 2t$ and $D_C^{1/2} t^2 = \frac{1}{\sqrt{\pi}} \int_0^t \frac{2s}{\sqrt{t-s}}\,ds = \frac{8\,t^{3/2}}{3\sqrt{\pi}}$.

**Constructive.** From Table 5, $c(2) = \frac{8}{3\sqrt{\pi}}$, hence $H_K(t^2) = \frac{8\,t^{3/2}}{3\sqrt{\pi}}$.

In all three ways $D_{RL}^{1/2} t^2 = D_C^{1/2} t^2 = H_K(t^2)$. It should be said here that, once the classical gamma formula (5.3) is known, the classical computation also takes one line. The value of the constructive approach does not lie in shortening the computation; it lies in the fact that the student builds this monomial rule algebraically.

### 5.3. An example of difference: $f(t) = 1$

Examples of agreement alone may create the false impression that the RL and Caputo operators are always identical. A contrasting example is therefore needed — a constant.

**Riemann–Liouville.** $J_{0+}^{1/2}1 = \frac{1}{\sqrt{\pi}}\int_0^t \frac{ds}{\sqrt{t-s}} = \frac{2\sqrt{t}}{\sqrt{\pi}}$, whence

$$D_{RL}^{1/2}1 = \frac{1}{\sqrt{\pi t}}. \qquad (5.5)$$

**Caputo.** $f' = 0$, therefore

$$D_C^{1/2}1 = 0. \qquad (5.6)$$

**Constructive.** With the normalization $c(0) = 1/\sqrt{\pi}$, $H_K(1) = \frac{1}{\sqrt{\pi t}}$.

Thus $D_{RL}^{1/2}1 = H_K(1) \neq D_C^{1/2}1$. The student sees that the names of the operators are not a formal distinction — the definitions really give different results. This is exactly where the choice of normalization in §4.5 becomes visible: $c_0 = 1/\sqrt{\pi}$ means reproducing the RL behaviour on constants, whereas the Caputo behaviour ($c_0 = 0$) contradicts the composition requirement (4.6).

The difference is explained by a single formula. If $f$ is absolutely continuous on $[0, T]$, then

$$D_C^{1/2}f = D_{RL}^{1/2}[f - f(0)] = D_{RL}^{1/2}f - \frac{f(0)}{\sqrt{\pi t}}. \qquad (5.7)$$

The two operators differ only by the contribution of the initial value: on monomials with $\beta > 0$ we have $f(0) = 0$ and they coincide (5.3), whereas on a constant $f(0) = 1$ and the difference is exactly $\frac{1}{\sqrt{\pi t}}$.

**Higher order.** The difference is not confined to order one-half. The Caputo derivative of order $\frac{3}{2}$ ($D_C^{3/2}f = J_{0+}^{1/2}f''$) vanishes on constant and linear functions: $D_C^{3/2}1 = D_C^{3/2}t = 0$, whereas the RL monomial rule gives $D_{RL}^{3/2}1 = -\frac{1}{2\sqrt{\pi}}\,t^{-3/2}$ and $D_{RL}^{3/2}t = \frac{1}{\sqrt{\pi t}}$. This shows the student the importance of the order, the initial conditions and the regularity, and that the choice of a classical operator must match the conditions of the problem.

### 5.4. Criteria for a meaningful comparison

The comparison of the three approaches should not be based on the length of the computation: as we saw in §5.2, with the known gamma formula the classical computation is as short as the constructive one. The paper therefore does not claim that the constructive method is faster or simpler. Its methodological purpose is that the process of construction itself becomes the object of learning. For the comparison we propose six criteria:

1. **Mathematical correctness** — does the result agree with the corresponding classical definition?
2. **Domain** — is the use of the operator on the given function legitimate?

3. **Route of computation** — is the integral definition, a ready-made gamma formula, or an algebraically constructed coefficient used?
4. **Conceptual understanding** — can the student explain the origin of the result?
5. **Independent check** — can the result be verified by another method?
6. **Instructional difficulty** — what prior knowledge does each route require, and where might help be needed?

These criteria underlie the problems of Section 6 and the rubric of Section 9. The constructive approach is regarded here as a possible initial route into fractional differentiation that complements, rather than replaces, the classical integral theory.

## 6. Applied problems: classical and constructive solutions side by side

This section contains seven problems (A–G). Each serves one *new* instructional aim not present in the preceding problems (Table 7). Computations already carried out in Sections 4 and 5 or in earlier problems are not repeated — we use the ready result with a reference. The problems are arranged in order of increasing conceptual difficulty: computation and checking → solving an equation → initial condition → modelling → general order → finite sums.

The constructive computations are restricted to $V_S$ (to $V_S^+$ where needed, §4.7); infinite series are not considered. A "classical solution" means a computation using the integral definition or formulas derived from it. As noted in §5.4, with the known gamma formula both routes may be equally short; the aim of the problems is not to compare speed. We use the values

$$\Gamma\left(\frac{1}{2}\right)=\sqrt{\pi},\quad \Gamma\left(\frac{3}{2}\right)=\frac{\sqrt{\pi}}{2},\quad \Gamma\left(\frac{5}{2}\right)=\frac{3\sqrt{\pi}}{4},\quad \Gamma\left(\frac{7}{2}\right)=\frac{15\sqrt{\pi}}{8},\quad \Gamma\left(\frac{9}{2}\right)=\frac{105\sqrt{\pi}}{16} \qquad (6.1)$$

and the monomial rule for the fractional integral

$$J_{0+}^{\alpha}t^{\beta}=\frac{\Gamma(\beta+1)}{\Gamma(\beta+\alpha+1)}\,t^{\beta+\alpha},\qquad \alpha>0,\ \beta>-1,\ t>0. \qquad (6.2)$$

### 6.1. Problem A — checking the composition

**Problem.** Verify that $H_K\big(H_K(t^2)\big)=D(t^2)$.

**New aim:** the difference between checking one example and the general justification obtained from the recurrence.

The half-derivative of $t^2$ has already been computed in three ways in §5.2: $H_K(t^2)=\frac{8}{3\sqrt{\pi}}\,t^{3/2}$. The second application, from Table 5, $c\left(\frac{3}{2}\right)=\frac{3\sqrt{\pi}}{4}$:

$$H_K\big(H_K(t^2)\big)=\frac{8}{3\sqrt{\pi}}\cdot\frac{3\sqrt{\pi}}{4}\,t=2t=D(t^2). \qquad (6.3)$$

**Domain.** $t^2$ and the intermediate $t^{3/2}$ both belong to $V_S$, so the composition is legitimate.

**Instructional role.** (6.3) is a special case of (4.5), $c(2)\,c\left(\frac{3}{2}\right) = 2$. The student is asked: does checking one monomial prove (4.1) on the whole of $V_S$? Answer — no; the general justification is provided by the recurrence and linearity. APOS: action → process.

### 6.2. Problem B — Abel's equation of the first kind

**Problem.** Find $u(t)$ if

$$J_{0+}^{1/2}u(t) = t^3, \qquad \text{i.e.} \qquad \frac{1}{\sqrt{\pi}}\int_0^t \frac{u(s)}{\sqrt{t-s}}\,ds = t^3. \qquad (6.4)$$

**New aim:** the operator is used for the first time to find an *unknown* function.

**Classical solution.** By the semigroup property $J_{0+}^{1/2}J_{0+}^{1/2} = J_{0+}^{1}$; applying $J_{0+}^{1/2}$ to both sides and differentiating gives $u = \frac{d}{dt}\left[J_{0+}^{1/2}t^3\right]$. By (6.2), $J_{0+}^{1/2}t^3 = \frac{\Gamma(4)}{\Gamma(9/2)}\,t^{7/2} = \frac{32}{35\sqrt{\pi}}\,t^{7/2}$, hence

$$u(t) = \frac{16}{5\sqrt{\pi}}\,t^{5/2}. \qquad (6.5)$$

**Constructive solution.** By §4.9, if $J_K u = g \in V_S$, then $u = H_K g$. From Table 5, $c(3) = \frac{16}{5\sqrt{\pi}}$, i.e. $u = H_K(t^3)$, and the same (6.5) is obtained in one step.

**Check.** $J_{0+}^{1/2}u = \frac{16}{5\sqrt{\pi}} \cdot \frac{\Gamma(7/2)}{\Gamma(4)}\,t^3 = t^3$.

**Domain.** $t^{5/2}$ is integrable near zero and the right-hand side lies in $V_S$. Checking the solution obtained must be distinguished from a general inversion theorem on a broad function space.

**Instructional role.** The third situation of the transition from process to object (§2.2): the operator appears in an equation. APOS: process → object.

### 6.3. Problem C — a Caputo initial value problem

**Problem.** Solve

$$D_C^{1/2}y(t) = t^2, \qquad y(0) = 0. \qquad (6.6)$$

**New aim:** the role of the initial condition and the link between the type of operator and the initial condition.

**Classical solution.** Under regularity conditions $J_{0+}^{1/2}\left(D_C^{1/2}y\right) = y(t) - y(0)$. Hence $y = J_{0+}^{1/2}t^2$, and by (5.4)

$$y(t) = \frac{16}{15\sqrt{\pi}}\,t^{5/2}. \qquad (6.7)$$

**Constructive solution.** We seek $y = At^\beta$. Exponent: $\beta - \frac{1}{2} = 2$, i.e. $\beta = \frac{5}{2}$. On monomials with positive exponent the Caputo operator coincides with $H_K$ (5.3), so from Table 5 $H_K(t^{5/2}) = \frac{15\sqrt{\pi}}{16} t^2$, whence $A = \frac{16}{15\sqrt{\pi}}$ — the same as (6.7). The initial condition holds, since $t^{5/2} \to 0$ as $t \to 0^+$.

**Check.** $y' = \frac{8}{3\sqrt{\pi}} t^{3/2}$ and $D_C^{1/2} y = J_{0+}^{1/2} y' = \frac{8}{3\sqrt{\pi}} \cdot \frac{\Gamma(5/2)}{\Gamma(3)} t^2 = t^2$.

**Role of the constant term.** If the solution is sought in the form $y = \kappa + At^{5/2}$, the Caputo operator annihilates $\kappa$: the equation does not determine $\kappa$; it is fixed by the initial condition ($y(0) = 0 \Rightarrow \kappa = 0$). With the RL operator the constant would additionally produce $\frac{\kappa}{\sqrt{\pi t}}$ (5.7). Thus the type of the operator and the form of the initial condition are linked.

**Domain.** $y$ is absolutely continuous on every $[0, T]$, so the Caputo definition applies.

**Instructional role.** The student constructs the unknown function, chooses the operator, and interprets the initial condition. APOS: object → schema.

### 6.4. Problem D — a normalized heat-flux model

**Context.** In some boundary value problems of heat conduction in a semi-infinite medium, the boundary temperature and the heat flux are related by fractional operators [14, 16]. In dimensionless variables we take the model

$$\theta(t) = J_{0+}^{1/2} p(t), \qquad (6.8)$$

where $\theta$ is the change of boundary temperature and $p$ is the normalized flux. (6.8) is taken as the starting model; its derivation from the heat equation is not the subject of this article.

**Problem.** Find $p(t)$ in two cases: (a) $\theta(t) = t$; (b) $\theta(t) = \theta_0$, $t > 0$ (temperature jump).

New aim: understanding the interplay between the integral model, the initial conditions, and the admissible class of functions when choosing a fractional operator.

(a) The mechanism of the solution is the same as in Problem B: $p = H_K\theta$; from Table 5, $p(t) = \frac{2\sqrt{t}}{\sqrt{\pi}}$; check: $J_{0+}^{1/2} p = \frac{2}{\sqrt{\pi}} \cdot \frac{\Gamma(3/2)}{\Gamma(2)} t = t$.

(b) With the RL and the constructed operator

$$p(t) = \frac{\theta_0}{\sqrt{\pi t}}. \qquad (6.9)$$

For the temperature jump, the flux obtained with the Riemann–Liouville operator is singular as $t \to 0^+$ but integrable near zero. The Caputo half-derivative gives zero on a constant function; the expression so obtained is not a solution of the given integral equation. The difference between the operators should therefore not be seen as an arbitrary choice of physical result: it is tied to the formulation of the model, the initial conditions, and the domain of the operator used.

Instructional role. The student checks why the Riemann–Liouville result agrees with the original integral model and why it cannot be replaced by the Caputo derivative of the constant without further justification. APOS: object → schema.

### 6.5. Problem E — the derivative of order $\frac{3}{2}$ and commuting operators

**Problem.** Find $D_{0+}^{3/2}t^3$ and check the result by compositions.

**New aim:** the transition from order one-half to a general order, and the conditions under which operators commute.

**Classical solution.** $D_{0+}^{3/2}t^3 = \frac{d^2}{dt^2}\left[J_{0+}^{1/2}t^3\right]$, and $J_{0+}^{1/2}t^3$ is already known from Problem B. Differentiating twice, or directly by (4.14),

$$D_{0+}^{3/2}t^3 = \frac{\Gamma(4)}{\Gamma(5/2)}\,t^{3/2} = \frac{8}{\sqrt{\pi}}\,t^{3/2}. \qquad (6.10)$$

**Constructive check in three ways.** $D\big(H_K(t^3)\big)$, $H_K\big(D(t^3)\big) = H_K(3t^2)$ and $H_K\left(H_K\big(H_K(t^3)\big)\right) = c(3)\,c\left(\frac{5}{2}\right)\,c(2)\,t^{3/2}$ — all three give $\frac{8}{\sqrt{\pi}}\,t^{3/2}$. All intermediate monomials lie in $V_S$, so the compositions are legitimate.

**Counterexample.** On a constant $D(H_K 1) = D\left(\frac{t^{-1/2}}{\sqrt{\pi}}\right) = -\frac{t^{-3/2}}{2\sqrt{\pi}} = D_{0+}^{3/2}1$, but $H_K(D1) = 0$. Thus $D \circ H_K \neq H_K \circ D$ on constants. This is the other side of the example in §4.8 ($D^{1/2}D \neq D^{3/2}$): commuting operators and adding orders are not admissible without checking the domains.

**Instructional role.** Computing the result is secondary here; the main point is understanding the conditions for composition and the counterexample. APOS: schema.

### 6.6. Problems F and G — finite sums

**New aim:** extending the operator constructed on monomials to finite sums by linearity. The mechanism of both problems is already known (B and C); only linearity and, in G, one additional observation are new.

**Problem F.** Find $u(t)$ if $J_{0+}^{1/2}u(t) = t + t^2 + t^3$. By linearity $u = H_K(t) + H_K(t^2) + H_K(t^3)$ (classically — $D_{0+}^{1/2}$ applied to each term), and from Table 5

$$u(t) = \frac{2}{\sqrt{\pi}}\,t^{1/2} + \frac{8}{3\sqrt{\pi}}\,t^{3/2} + \frac{16}{5\sqrt{\pi}}\,t^{5/2}. \qquad (6.11)$$

A term-by-term check with $J_{0+}^{1/2}$ gives $t + t^2 + t^3$. Only a finite combination is used; no conclusion about infinite series follows from this.

**Problem G.** Solve $D_C^{1/2}y = 1 + t + t^2$, $y(0) = 0$. Classically $y = J_{0+}^{1/2}(1 + t + t^2)$. Constructively we seek $y = At^{1/2} + Bt^{3/2} + Ct^{5/2}$; from Table 5 $A\frac{\sqrt{\pi}}{2} = 1$, $B\frac{3\sqrt{\pi}}{4} = 1$, $C\frac{15\sqrt{\pi}}{16} = 1$, that is,

$$y(t) = \frac{2}{\sqrt{\pi}} t^{1/2} + \frac{4}{3\sqrt{\pi}} t^{3/2} + \frac{16}{15\sqrt{\pi}} t^{5/2}. \qquad (6.12)$$

$y(0) = 0$ holds; near zero $y'$ contains an integrable term proportional to $t^{-1/2}$, so $y$ is absolutely continuous and the Caputo definition is legitimate.

**Additional observation from G.** The constant term 1 on the right-hand side arises from the half-derivative of $t^{1/2}$, and *not* from including a constant in the solution — the Caputo operator annihilates constants (Problem C). This is a source of a typical student error (Section 7).

**Instructional role.** F — understanding linearity (APOS: object); G — combining the ideas of B, C and F (APOS: schema).

### 6.7. Overview of the problem set

*Table 7. New instructional aims of the applied problems*

| Problem | Mathematical action | New instructional aim | Ready result used | Expected APOS emphasis |
|---|---|---|---|---|
| A | checking the composition | example ≠ general justification | §5.2, Table 5 | action → process |
| B | Abel equation | the operator used to find an unknown function | §4.9 | process → object |
| C | Caputo initial value problem | role of the initial condition | (5.4), (5.7) | object → schema |
| D | heat-flux model | physical consequence of the choice of operator | B, §5.3 | object → schema |
| E | derivative of order $\frac{3}{2}$ | general order, conditions for commuting | B, (4.14), §4.8 | schema |
| F | Abel equation with a sum | linearity | B | object |
| G | Caputo problem with a sum | combining ideas; origin of the constant term | C, F | schema |

This sequence starts instruction with a simple action and leads to situations in which the student must simultaneously take into account domains, inversion, composition, and initial conditions. The indicated APOS stages are tentative reference points of the instructional design, not empirically established levels of students' mental development.

## 7. Proposed classroom implementation: anticipated errors and the lecturer's response

### 7.1. Instructional blocks

For classroom work the seven phases of Table 4 can be grouped into four blocks:

- **Block I — formulating the problem** (phases 1–2, §4.1–4.2). The student recalls differentiation and the repeated application of an operator, formulates the requirement $H \circ H = D$, chooses the monomial class, and arrives at the form of the operator (4.3) through explicitly declared assumptions.
- **Block II — recurrence and families** (phases 3–4, §4.3–4.4). The student applies the operator twice, derives (4.5), computes specific coefficients, and discovers the free parameter.
- **Block III — normalization and comparison** (phase 5, §4.5–4.6, Section 5). The student understands the need for an additional condition, chooses the normalization compatible with RL, and compares the three operators.
- **Block IV — domains and problems** (phases 6–7, §4.7–4.9, Section 6). The student checks the domains and applies the construction in problems.

A block does not necessarily mean a single session: its length depends on the students' prior knowledge, preparation, and need for help. Phase 7 (general orders, the half-integral) may be left as independent work for interested students. This organization is a methodological design; its optimal duration and effectiveness have not been established empirically.

### 7.2. Anticipated errors and their didactic use

In constructive teaching an error is not merely a wrong result; it points to a cognitive obstacle that arose at a specific stage of building knowledge. The lecturer's aim is therefore to identify the cause of the error and to ask a question that allows the student to check their own answer. Nine anticipated errors are briefly described below; a summary is given in Table 8. Each is linked to a phase of Table 4 or to a specific place in Sections 5–6.

**1. Taking the requirement as a statement of existence (phase 1).** The student writes $H \circ H = D$ as a known fact. Lecturer's question: "Do we know that such an $H$ exists? What must we specify in order to find out?"

**2. Taking an assumption for a derived result (phase 2).** The student thinks that the monomial form and linearity follow from (4.1). Lecturer's question: "Did we obtain this from $H \circ H = D$, or did we choose it? Is another choice possible?"

**3. Mechanically halving the coefficient (phase 3).** The student proposes $H(x^\beta) = \frac{\beta}{2} x^{\beta-1/2}$. The change of exponent is correct, the coefficient is not. The lecturer suggests applying the rule twice to $x$: $H(x) = \frac{1}{2} x^{1/2}$, $H(H(x)) = \frac{1}{2} \cdot \frac{1}{4} = \frac{1}{8}$, whereas $D(x) = 1$. The student sees that the guess does not satisfy the composition. Question: "What must change in the coefficient so that applying the rule twice gives the derivative?"

**4. Taking $c_0$ to be uniquely determined (phases 4–5).** The student thinks that (4.5) determines all coefficients uniquely, or that $c_0 = 1/\sqrt{\pi}$ follows from the recurrence. The lecturer suggests two different values of $c_0$ and the computation of several coefficients: the coefficients change, while (4.5) still holds. Question: "What additional condition do we need in order to choose a specific operator?"

**5. "Every derivative sends a constant to zero" (Section 5).** By analogy with ordinary differentiation. The lecturer suggests comparing (5.5) and (5.6): $D_{RL}^{1/2}1 = \frac{1}{\sqrt{\pi t}}$, $D_{C}^{1/2}1 = 0$. Question: "Which definition gives zero, and why is this not common to all fractional derivatives?"

**6. Ignoring the domain of the composition (phases 6–7).** The student applies the operator twice without checking where the intermediate function belongs. The lecturer returns to §4.7: $H_K\left(H_K\left(x^{-1/2}\right)\right) = 0$, whereas $D\left(x^{-1/2}\right) = -\frac{1}{2}x^{-3/2}$. Question: "Are both operations defined, and on which space have we established $H_K^2 = D$?"

**7. Taking computational simplicity for a universal advantage (§5.4).** The student concludes that the constructive route is always simpler. The lecturer compares two situations: a classical computation from the integral definition and a classical computation with the known gamma formula. Question: "If the gamma formula is already known, how different are the numbers of steps?" The correct conclusion is qualified: the constructive route may have an advantage in understanding the origin of the formula, but not a universal computational advantage.

**8. Treating negative exponents as automatically inadmissible (§5.1).** A negative exponent must be distinguished from integrability at zero: for the RL integral $\beta > -1$ is sufficient, whereas the standard Caputo definition also requires integrability of $f'$. Question: "Why may $t^{-1/2}$ belong to the relevant RL class and yet fail to meet the Caputo requirements?"

**9. Using the table outside the monomial class (phase 6).** The student applies the coefficients of Table 5 to a function that does not belong to $V_S$, e.g. to $\sin t$. Passing to such functions requires a separate theoretical extension (infinite series, convergence, analysis of the domain). Lecturer's question before any computation: "Does the function belong to the class on which our operator is defined?"

**Addendum from Problem G.** In a Caputo problem, because of the constant term on the right-hand side, the student may include a constant in the solution (§6.6). This is the reverse form of error 5 and is revealed by the same question: "What does the Caputo operator do to a constant, and where, then, must 1 come from?"

*Table 8. Anticipated student errors and the lecturer's proposed response*

| # | Symptom of the error | Lecturer's question | What the student should show |
|---|---|---|---|
| 1 | $H \circ H = D$ is taken as a statement of existence | Do we know that such an $H$ exists? | separating a requirement from a proof |
| 2 | the monomial form is taken as derived | Did we obtain this or choose it? | naming the assumption explicitly |
| 3 | coefficient $\beta/2$ | Apply the rule twice to $x$ — do you get 1? | checking the composition, obtaining (4.5) |
| 4 | $c_0$ is taken as uniquely determined | What happens for two values of $c_0$? | explaining the free parameter and the need for normalization |

| # | Symptom of the error | Lecturer's question | What the student should show |
|---|---|---|---|
| 5 | RL and Caputo are taken to be identical; a constant in the solution (G) | What do we get on a constant? | correct comparison of the two definitions, (5.7) |
| 6 | composition is used unconditionally | To which space does the intermediate function belong? | checking the domains |
| 7 | the constructive route is taken to be always simpler | What changes if the gamma formula is known? | taking the conditions of comparison into account |
| 8 | a negative exponent is treated as automatically inadmissible | What happens when $-1 < \beta < 0$? | justifying the integrability conditions |
| 9 | the table is applied to any function | Does the function belong to $V_S$? | stating the domain before computing |

### 7.3. From discussing an error to an instructional decision

The aim of independent work. The student should not only arrive at the correct answer but also be able to state his or her own assumption, check an algebraic step, explain the cause of a wrong conclusion, and correct the error independently. A targeted hint is used only when a specific obstacle has been identified. Improvement of these skills is a hypothesis to be tested in future empirical research, not a result established in this paper.

The questions of Table 8 are not an empirically validated diagnostic instrument. They are based on the structure of the construction and on presumed cognitive difficulties, so that the lecturer has a prepared response when a student cannot take a step or obtains a wrong result. Identifying an error is not the end of the instructional action: it becomes the basis for a further question, a move to a simpler task, or a renewed activation of prior knowledge (§2.1, Vygotsky).

It is precisely this principle that connects Part I with Part II. In Part I these questions are methodological tools of the proposed classroom implementation; the questions and the levels of help have already been implemented in the program "Sokrate AI — Kapanadze" (§1.5; [6]), which identifies the type of the student's error, gives graded hints, and records progress; Part II is devoted to its description. The general principles for the role of artificial intelligence in this part are briefly set out in Section 8.

## 8. The role of artificial intelligence: principles and limits

This section is deliberately short. It sets out only those principles that bear directly on the instructional sequence of Part I. These principles underlie the program "Sokrate AI — Kapanadze" (§1.5; [6]); the scenario of the Socratic dialogue, its samples, the levels of help, and the program's diagnostics are described in Part II of the study.

### 8.1. The function of AI and the roles of the three participants

AI is considered here as an instructional assistant, not as a substitute for mathematical reasoning; the discussion concerns methodological possibilities, not the effectiveness of a finished digital tutor. The roles are distributed as follows: the **lecturer** defines the instructional aim, the logic of the questions (Tables 4 and 8), the admissible help, and the criteria of correctness; the **student** carries out the main cognitive work; **AI** may help with guiding questions and by checking intermediate steps — only within the limits designed by the lecturer.

The educational opportunities and risks of large language models are widely discussed [20]; specialized models for Socratic mathematics teaching already exist [4, 5]. Empirical studies show that the outcome depends not only on access to AI but also on the pedagogical organization of the help it provides.

In the randomized crossover experiment of Kestin and co-authors [21] (an introductory physics course at Harvard University, autumn 2023, 194 students), students learned more in a lesson taught by an AI tutor with a research-based design than in active in-class instruction. In the field experiment of Bastani and co-authors [22] (about a thousand high-school students, mathematics), unrestricted access to GPT-4 improved performance during practice, but after access was removed the students performed worse than the control group; a version with learning safeguards largely eliminated this negative effect. Both studies concern a different subject and context; **neither tested the half-derivative design proposed here**. They only justify why the restrictions described below are needed.

### 8.2. Two modes

Methodologically we distinguish two modes.

- **Discovery mode** (phases 1–5). AI does not give the student the recurrence, the value of $c_0$, or the gamma formula before the student has taken the corresponding step. Its function is to ask the questions of Tables 4 and 8 and, when needed, to move to a simpler question (§2.1).
- **Verification mode** (phase 6 and Section 6). Checking a formula, a table, or computations that have already been derived is permitted — subject to the safeguards described below.

This is a proposed pedagogical requirement, not a software feature tested in this article.

### 8.3. Safeguards against computational errors

Every computation by AI requires independent checking: a language model may confuse $\Gamma\left(\frac{5}{2}\right)$ with $\Gamma\left(\frac{7}{2}\right)$, lose a factor $\sqrt{\pi}$, or attribute the RL value to the Caputo derivative of a constant. The topic depends on exact values of the gamma function, so we propose the following rules; none of them relies on trust in AI.

1. **Reference list.** The values (6.1) are supplied by the lecturer in advance; every gamma value obtained in the dialogue is compared with this list by the student, not by AI.
2. **Derivation instead of recall.** In verification mode AI obtains every half-integer value from $\Gamma\left(\frac{1}{2}\right) = \sqrt{\pi}$, writing out the steps of $\Gamma(z+1) = z\Gamma(z)$, so that an error is visible in a specific line.

3. **Composition test.** Every pair of neighbouring coefficients must satisfy (4.5). This is an internal check of the construction, done by hand in one line, and it immediately reveals a wrong coefficient.
4. **Inversion test.** Every solution of an equation of the form $J_{0+}^{1/2}u = g$ is substituted back using (6.2), as in Problems B and D.
5. **Constant test.** Before a problem involving the Caputo operator the question is asked: "What is the half-derivative of 1 for each of the three operators?" An answer that does not distinguish Caputo from RL is a signal to stop and check again.
6. **Numerical confirmation.** Where a computer algebra system is available, a single numerical evaluation of the integral (5.1) at one point is compared with the closed form.
7. **Responsibility of the lecturer.** Final correctness is confirmed by the lecturer; a proof generated by AI is not accepted as justification of a step in assessed work.

These rules also have an intended didactic effect: the student practises checking as an independent mathematical activity — which is one of the criteria in Section 9.

### 8.4. Differentiated support and ethics

For a student who needs more support, AI can offer separately worked steps, e.g. the computation of $c\left(\frac{1}{2}\right)$, $c(1)$ and $c\left(\frac{3}{2}\right)$; for a strong student, questions about the freedom of normalization, composition of order $\frac{3}{2}$, or the case $-1 < \beta < 0$ (phase 7). The level of support should be based on the student's actual answers, not on unverified assumptions about the student.

The lecturer should request only the minimum necessary data, should not upload information that identifies the student, and should clearly explain the role of AI in assessment [23]. For further research it is of interest how far AI used according to these principles helps the student to discover errors independently, to justify reasoning, and to solve problems without help. Separate indicators of knowledge components (declarative, procedural, metacognitive) and the history of hint use are discussed in the design of Part II; they are not yet validated measures.

## 9. Assessment framework and further empirical research

### 9.1. Assessment criteria

From the cognitive-constructivist point of view it matters not only whether the student obtained the correct answer but also how he or she formulated the problem, whether the relationship was derived independently, and how thoughtfully the result was checked. We therefore propose a rubric for mathematical reasoning, each criterion of which is linked to a specific phase of the construction (Table 4) and to an anticipated error (Table 8). A three-level scale is used:

- **0** — the student cannot carry out the action or cannot give a correct explanation;
- **1** — carries it out partially or with targeted help;
- **2** — carries it out independently, correctly, and with justification.

This scale separates a result obtained independently from an action performed with help.

*Table 9. Assessment rubric for the constructive teaching of the half-derivative (proposed, not validated)*

| **Criterion** | **Phase / error** | **0** | **1** | **2** |
|---|---|---|---|---|
| K1. Formulating the problem | 1 / #1 | cannot formulate $H \circ H = D$ | formulates it with a hint | formulates it independently and explains that it is a requirement, not a proof |
| K2. Separating assumption from conclusion | 2, 5 / #2, #4 | does not distinguish what was chosen from what was derived | distinguishes them in some cases | names the status of the formula at every step (§3.1) |
| K3. Recurrence relation | 3 / #3 | does not obtain the correct relation | obtains it with errors or with help | derives (4.5) independently and states the condition $\beta \geq \frac{1}{2}$ |
| K4. Understanding the normalization | 4–5 / #4 | does not notice the freedom of $c_0$ | acknowledges the need for an additional condition but cannot explain why | justifies the non-uniqueness and explains the choice of normalization |
| K5. Domain | 6–7 / #6, #8, #9 | ignores the restrictions | notes some restrictions but does not fully check intermediate functions | independently checks the domains and the integrability conditions |
| K6. Comparison with the classical operators | Section 5 / #5, #7 | does not distinguish RL from Caputo | notices the difference in a specific example | justifies the conditions of agreement and difference from the definitions |
| K7. Checking the result | Section 6; §8.3 | does not check the answer | the check is incomplete or done with help | independently checks the result in the original problem or by another method |

### 9.2. A criterion profile instead of a total score

Until the rubric has been validated, the result should be reported as a **criterion profile** — a sequence of seven scores $(K1, \dots, K7)$ — and not as a sum or a percentage. The reason is that the same sum can hide entirely different difficulties. For example, two students may have the profiles

$$A = (2,\ 2,\ 2,\ 0,\ 1,\ 2,\ 1), \qquad B = (2,\ 1,\ 1,\ 2,\ 2,\ 1,\ 1).$$

Both have a sum of 10, but A computes the coefficients correctly and cannot explain the meaning of the normalization, whereas B understands the operator construction well but makes algebraic errors. The lecturer's next step differs in the two cases: for A — question 4 of Table 8; for B — the composition test (§8.3). A total score loses this information. Moreover, a score of 1 combines two different situations (partial performance and performance with help); when recording the profile it is advisable to indicate which of them applies. A sum or a percentage may be used only once the unidimensionality of the rubric has been established empirically.

### 9.3. Using the rubric in instruction

The rubric is used after instruction to assess understanding of the construction. The initial diagnostic is a procedure with a different purpose: it should cover only the knowledge expected before instruction (§2.1) — ordinary differentiation, the power rule, composition of operators, and algebraic manipulation. It is not admissible to require knowledge of a formula the student has not yet studied and to count its absence as a gap in prior knowledge.

### 9.4. How the model could be tested

Empirical testing of the model requires two stages: first the validation of the instruments themselves, then the study of the instructional effect.

**Validation of the instruments.** (a) *Content validity*: expert review of the rubric and of the questions of Table 8 by specialists in mathematics and mathematics education. (b) *Inter-rater agreement*: two independent raters score the same written work, and agreement is computed separately for each criterion. Since 0–1–2 is an ordinal scale, instead of the ordinary Cohen's kappa one uses the **quadratically weighted kappa** (or the ordinal version of Krippendorff's $\alpha$): these indices take into account that a disagreement between 0 and 2 is more serious than one between 0 and 1. (c) *Pilot implementation* with a small group, after which unclear level descriptors are refined.

**Data sources.** A preliminary diagnostic (prior knowledge); post-instruction written tasks scored with the rubric; analysis of written reasoning — by types of arguments and errors (Table 8), with the students' consent; conceptual interviews with selected students — to find out whether they treat the half-derivative as an operator in its own right (the indicators of §2.2); a delayed assessment — to check retention and transfer to new problems.

**Comparative design.** A quasi-experimental comparison is possible: one group studies with the constructive sequence, the other with the definition-first presentation (Table 3). Such a comparison requires equating the groups' prior knowledge, instructional time and assessment instruments, and, when different lecturers are involved, taking into account the influence of teaching style.

**Hypotheses to be tested.**

- **H1.** Students prepared with the constructive sequence justify the origin of the coefficients (K3, K4) better than those who learn to apply a ready-made formula.

- **H2.** Instruction organized through guiding questions supports students in detecting and correcting their own errors (K7).
- **H3.** Treating the half-derivative as an operator in its own right improves the ability to check the conditions for composition and the domains (K5).
- **H4.** Explicit labelling of the logical status of formulas improves the separation of assumption from conclusion (K2).

These are hypotheses for future research, not results of this article.

### 9.5. Connection with Part II

The rubric of Part I is tied to the specific mathematical goals of the half-derivative. Part II discusses a broader process-oriented diagnostic that the program "Sokrate AI — Kapanadze" [6] produces at the end of a session: a record of the help used, identification of independent actions, and indicators of declarative, procedural and metacognitive knowledge. From this point of view the rubric of Part I is the initial instrument, and the diagnostics of Part II are its possible extension. In both cases the reliability and validity of the criteria must be tested independently; the present paper does not claim that these instruments are validated or that they establish the effectiveness of the model.

## 10. Discussion, limitations and conclusion

### 10.1. Discussion of the results

The paper turns the author's algebraic-operator construction [1, 2] into an instructional sequence in which the student moves step by step from prior knowledge to the monomial formula for the half-derivative. The mathematical and the didactic contributions are kept apart (Table 2): the mathematical result — the construction of the operator on $V_S$, the recurrence, the free normalization, and the correspondence with the classical operators — is taken from earlier work; the contribution of this article is its organization into phases, diagnostic questions, problems, and assessment criteria. The separation matters because the correctness of a mathematical result does not by itself mean that a particular way of teaching it is effective.

The central moment of the sequence is normalization. Within the chosen assumptions the composition requirement yields the recurrence but does not determine a unique operator: the student faces a situation where the derived result is correct, yet an additional choice is needed for the final operator. This is exactly where the difference between a mathematical conclusion and a convention becomes visible — a difference that the labelling of the logical status of formulas (§3.1) maintains throughout the text.

The second main instructional moment is the domain: the operator is defined on $V_S^+$, whereas the identity $H_K \circ H_K = D$ holds only on $V_S$ (§1.2, §4.7). The comparison with the classical operators (Section 5) showed that the result of an operator depends not only on the notation for the order but also on the definition, the class of functions, and the initial value (5.7): the normalized $H_K$ agrees with RL for all exponents of $S$, while Caputo agrees with them for positive exponents and differs on constants.

The problems of Section 6 are built so that each has one new instructional aim and the student gradually moves from computation to treating the operator as an object in its own right. This transition has not been confirmed with real students. The comparison of the classical and constructive routes shows that on monomials and finite sums integral computations can be replaced by algebra, but with the known gamma formula this difference shrinks. The paper therefore does not claim a universal computational advantage for the constructive route; its purpose is to make the origin of the formula and the conditions of its use understood.

### 10.2. Limitations

1. **Mathematical scope.** The construction is defined on $V_S$; on the extended space $V_S^+$ the composition identity no longer holds.
2. **Broader classes of functions.** Infinite series and general functions are not considered; such an extension requires a separate study of convergence, domains, and term-by-term application. Direct transfer of the monomial formulas to arbitrary functions is not admissible.
3. **Domains of the classical operators.** The RL and Caputo operators have different regularity requirements; the comparison (Section 5) is restricted to monomials with $\beta \geq 0$.
4. **Absence of empirical research.** No systematic experiment was carried out; cognitive load, instructional effectiveness, and the development of independent reasoning were not measured. The informal use of the program at SEU and conversations with colleagues (§1.5) are not treated as data; the student's question (§1.1) is a motivating episode.
5. **Indicators and instruments.** The APOS stages, the error catalogue, and the rubric were developed within the instructional design and are not validated instruments; results should therefore be reported as criterion profiles (§9.2).
6. **Artificial intelligence.** In Part I the role of AI is described only at the level of principles (Section 8); the program "Sokrate AI — Kapanadze" [6], its dialogues and diagnostics are described in Part II; the existence of the program does not by itself confirm the effectiveness of the model.

### 10.3. Achievement of the study objectives

- **O1 (instructional sequence)** is achieved at the level of design: the mathematical route of [1, 2] is organized into seven phases (Table 4), each with prior knowledge, an anticipated obstacle and a diagnostic question, and the formulas are labelled by logical status.
- **O2 (comparison with the classical operators)** is achieved within the class and conditions considered: the three operators are compared on the same set of monomials (Table 6), and their agreement, the difference on constants, and its explanation (5.7) are identified.
- **O3 (applied problems)** is achieved: seven problems with classical and constructive solutions, each with one new instructional aim (Table 7).
- **O4 (classroom implementation and assessment)** is achieved at the conceptual-methodological level: instructional blocks and nine anticipated errors (Table 8), principles and safeguards for the use of AI (Section 8), a seven-criterion rubric, and a testing plan (Table 9, §9.4).

Fulfilling these objectives does not mean that the practical effectiveness of the model has been established; that is the next stage of the research.

### 10.4. Conclusion

The paper proposes a design for teaching the half-derivative in which a student's question — whether a derivative of order one-half makes sense, raised independently and without any knowledge of fractional calculus — becomes the requirement $H \circ H = D$, then, through explicitly declared assumptions, a monomial construction, a recurrence relation, the discovery of a free normalization and, through the convention of compatibility with Riemann–Liouville, the gamma formula. Along this route the student sees the final formula not as a ready-made rule but as the result of interconnected steps, and learns what has been derived, what has been chosen, and where the resulting identity holds. The comparison with the classical operators and the system of problems complement this understanding, and the constructive route is proposed not as a replacement for the classical integral theory but as an introduction to it.

The model, the diagnostic questions, and the rubric are conceptual-methodological results; their effectiveness and reliability require empirical testing. This lecturer-designed sequence has already been implemented in the program "Sokrate AI — Kapanadze" [6], which helps the student answer the questions with graded hints and finally produces a diagnostic of their progress. Part II of the study is devoted to this program, its diagnostic logic, and the assessment plan.

### Acknowledgements

The author thanks fellow mathematicians at universities in Tbilisi and Telavi for friendly and critical conversations about whether this line of work is worth pursuing, and the students of SEU for the question with which this work began. These conversations do not constitute an assessment of the content or results of the paper; responsibility for any remaining shortcomings rests solely with the author.